
\magnification 1250
\font\san=cmssdc10
\def\sy{\hbox{\san\char 83}} 

\input amssym.def
\mathcode`A="7041 \mathcode`B="7042 \mathcode`C="7043
\mathcode`D="7044 \mathcode`E="7045 \mathcode`F="7046
\mathcode`G="7047 \mathcode`H="7048 \mathcode`I="7049
\mathcode`J="704A \mathcode`K="704B \mathcode`L="704C
\mathcode`M="704D \mathcode`N="704E \mathcode`O="704F
\mathcode`P="7050 \mathcode`Q="7051 \mathcode`R="7052
\mathcode`S="7053 \mathcode`T="7054 \mathcode`U="7055
\mathcode`V="7056 \mathcode`W="7057 \mathcode`X="7058
\mathcode`Y="7059 \mathcode`Z="705A
\def\spacedmath#1{\def\packedmath##1${\bgroup\mathsurround=0pt##1\egroup$}
\mathsurround#1 \everymath={\packedmath}\everydisplay={\mathsurround=0pt}}
 \spacedmath{2pt}

\def\note#1#2{\footnote{\parindent
0.4cm$^#1$}{\vtop{\sevenrm\baselineskip12pt\spacedmath{1pt}
\hsize15.5truecm\noindent
#2}}\parindent 0cm}
\def\up#1{\raise 1ex\hbox{\smallf@nt#1}}
\def\tx{\kern-1.5pt -}
\def\cqfd{\kern 2truemm\unskip\penalty 500\vrule height
4pt depth 0pt width 4pt\medbreak}  \def\virg{\raise .4ex\hbox{,}}
\def\decale#1{\smallbreak\hskip 28pt\llap{#1}\kern 5pt}
\def\no{n\up{o}\kern 2pt}
\def\ind{\par\hskip 1truecm\relax}

\parindent=0cm
\def\moins{\mathrel{\hbox{\vrule height 3pt depth -2pt
width 6pt}}}
\def\sdir_#1{\mathrel{
\mathop{\kern 0pt\oplus}\limits_{#1}}}

\def\rond{\kern 1pt{\scriptstyle\circ}\kern 1pt}

\def\iso{\mathrel{\mathop{\kern 0pt\longrightarrow }\limits^{\sim}}}

\def\Tr{\mathop{\rm Tr}\nolimits}

\vsize = 25truecm
\hsize = 16truecm
\hoffset = -.15truecm
\voffset = -.5truecm
\baselineskip15pt
\overfullrule=0pt
\centerline{\bf A Calabi-Yau threefold with non-Abelian fundamental
group} \smallskip
\centerline{Arnaud Beauville}\vskip1truecm
\ind The aim of this note is to answer a question of I.
Dolgachev by constructing a Calabi-Yau threefold whose fundamental group is the
quaternionic group $H$ with $8$ elements.  The construction is very reminiscent
of
Reid's unpublished construction of  a surface  with $p_g=0$, $K^2=2$ and
$\pi_1=H$; I explain below the link between the two problems.
\vskip1truecm
{\bf 1. The example}
\smallskip \ind Let $H=\{\pm1,\pm i,\pm j,\pm k\}$ be the quaternionic group,
and $V$  its regular representation. We denote by $\widehat{H}$ the
 group of characters $\chi:H\rightarrow {\bf C}^*$; it is isomorphic to
${\bf Z}_2\times{\bf Z}_2$. The group $H$ acts on ${\bf P}(V)\!$\note{1}{I use
Grothendieck's notation, i.e. $\scriptstyle{\bf P}(V)$ is the space of
hyperplanes
in $\scriptstyle V$.}\  and on $\sy^2V$; for each $\chi\in\widehat{H}$, we
denote by
$(\sy^2V)_\chi$ the eigensubspace of $\sy^2V$ with respect to $\chi$, i.e. the
space
of quadratic forms $Q$ on ${\bf P}(V)$ such that $h\cdot Q =\chi(h)\,Q$ for all
$h\in
H$. \smallskip  {\bf Theorem}$.-$ {\it For each $\chi\in\widehat{H}$, let
$Q_\chi$ be
a general element of $(\sy^2V)_\chi$.  The subvariety $\widetilde X$ of ${\bf
P}(V)$
defined by the $4$ equations $Q_\chi=0$ $(\chi\in\widehat{H})$ is a smooth
threefold,
on which the group $H$ acts freely. The quotient $X:=\widetilde X/H$ is a
Calabi-Yau
threefold with $\pi_1(X)=H$.} \ind Let me observe first that the last assertion
is an
immediate consequence of the others. Indeed, since $\widetilde X$ is a
Calabi-Yau
threefold, one has $h^{1,0}(\widetilde X)=h^{2,0}(\widetilde X)=$ $\chi({\cal
O}_{\widetilde X})=0$, hence  $h^{1,0}(X)=h^{2,0}(X)=\chi({\cal O}_{X})=0$.
This
implies $h^{1,0}(X)=1$, so there exists a nonzero holomorphic $3$\tx form
$\omega$ on
$X$;  since its pull-back to $\widetilde X$ is everywhere nonzero, $\omega$ has
the
same property, hence $X$ is a Calabi-Yau threefold. Finally $\widetilde X$ is a
complete intersection in ${\bf P}(V)$, hence simply connected by Lefschetz'
theorem,
so the fundamental group of $X$ is isomorphic to $H$.

 \ind So the problem is to prove that $H$ acts freely and $\widetilde X$ is
smooth. We
will need to write down explicitely the elements of  $(\sy^2V)_\chi$. As a
$H$\tx
module, $V$ is the direct sum of the $4$ one-dimensional representations of $H$
and
twice the irreducible two-dimensional representation $\rho$. Thus there exists
a
system of homogeneous coordinates $(X_1,X_\alpha,X_\beta,X_\gamma;Y,Z;Y',Z')$
such
that  $$\hss
g\cdot(X_1,X_\alpha,X_\beta,X_\gamma;Y,Z;Y',Z')
=(X_1,\alpha(g)X_\alpha,\beta(g)X_\beta,
\gamma(g)X_\gamma;\rho(g)(Y,Z);\rho(g)(Y',Z'))\ .\hss$$
\ind To be more precise, I denote by $\alpha$ (resp. $\beta$, resp. $\gamma$)
the
nontrivial character which is $+1$ on $i$ (resp. $j$, resp. $k$), and I take
for
$\rho$ the usual representation via Pauli matrices:
$$\hss\rho(i)\,(Y,Z)=(\sqrt{-1}\,Y,-\sqrt{-1}\,Z)\quad
\rho(j)\,(Y,Z)=(-Z,Y)\quad
\rho(k)(Y,Z)=(-\sqrt{-1}\,Z,-\sqrt{-1}\,Y)\ .\hss$$
 Then the general element $Q_\chi$ of $(\sy^2V)_\chi$ can be written
$$\eqalign{Q_1&=t^{1}_1\,X^2_1+t^{1}_2\,X^2_\alpha+t^{1}_3\,X^2_\beta+t^{1}_4\,
X^2_\gamma+t^{1}_5\,(
YZ'-Y'Z)\cr
Q_\alpha&=t^{\alpha}_1\, X_1X_\alpha +t^{\alpha}_2 \,X_\beta X_\gamma+
t^{\alpha}_3\,YZ+t^{\alpha}_4\,Y'Z'+t^{\alpha}_5\,(YZ'+ZY')\cr
Q_\beta&=t^{\beta}_1\, X_1X_\beta +t^{\beta}_2 \,X_\alpha X_\gamma+
t^{\beta}_3\,(Y^2+Z^2)+t^{\beta}_4\,(Y'^2+Z'^2)+t^{\beta}_5\,(YY'+ZZ')\cr
Q_\gamma&=t^{\gamma}_1\, X_1X_\gamma +t^{\gamma}_2\, X_\alpha X_\beta+
t^{\gamma}_3\,(Y^2-Z^2)+t^{\gamma}_4\,(Y'^2-Z'^2)+t^{\gamma}_5\,(YY'-ZZ')\
.}$$

  \ind For ${\bf t}:=(t^\chi_i)$ fixed, let  ${\cal X}_{\bf t}$ be the
subvariety of ${\bf P}(V)$ defined by the equations $Q_\chi=0$. Let us check
first that
the action of $H$ on ${\cal X}_{\bf t}$ is fixed point free for ${\bf t}$
general
enough. Since a point fixed by an element $h$ of $H$ is also fixed by $h^2$, it
is
sufficient to check that the  element  $-1\in H$ acts without fixed
point, i.e. that ${\cal X}_{\bf t}$ does not meet the linear spaces $L_+$ and
$L_-$
defined by  $Y=Z=Y'=Z'=0$ and $X_1=X_\alpha=X_\beta=X_\gamma=0$ respectively.
\ind Let $x=(0,0,0,0;Y,Z;Y',Z')\in {\cal X}_{\bf t}\cap L_- $.
 One of the
coordinates, say $Z$,  is nonzero; since $Q_1(x)=0$, there exists $k\in{\bf
C}$ such that $Y'=kY$, $Z'=kZ$.  Substituting in the equations
$Q_\alpha(x)=Q_\beta(x)=Q_\gamma(x)=0$ gives
$$(t^\alpha_3+t^\alpha_5k+t^\alpha_4k^2)\,YZ=(t^\beta_3+t^\beta_5k
+t^\beta_4k^2)
\,(Y^2+Z^2)= (t^\alpha_3+t^\alpha_5k+t^\alpha_4k^2)\,(Y^2-Z^2)=0$$ which has no
nonzero solutions for a generic choice of ${\bf t}$.
\ind Now let $x=(X_1,X_\alpha,X_\beta,X_\gamma;0,0;0,0)\in {\cal X}_{\bf t}\cap
L_+ $. As soon as the $t_i^\chi$'s are nonzero, two of the $X$\tx coordinates
cannot
vanish, otherwise all the coordinates would be zero. Expressing that
$Q_\beta=Q_\gamma=0$ has a nontrivial solution in $(X_\beta,X_\gamma)$ gives
$X_\alpha^2$ as a multiple of $X_1^2$, and similarly
for $X_\beta^2$ and $X_\gamma^2$. But then  $Q_1(x)=0$ is impossible for a
general
choice of ${\bf t}$.
\ind Now we want to prove that ${\cal X}_{\bf t}$ is smooth for ${\bf t}$
general
enough.
 Let ${\cal
Q}=$ $\sdir_{\chi\in\widehat{H}}(\sy^2V)_\chi$; then
${\bf t}:=(t^\chi_i)$ is a system of
coordinates on ${\cal Q}$. The  equations $Q_\chi=0$ define a subvariety ${\cal
X}$ in
${\cal Q}\times{\bf P}(V)$, whose fibre  above a point
${\bf t}\in{\cal Q}$ is ${\cal X}_{\bf t}$.
Consider the second projection $p:{\cal X}\rightarrow {\bf P}(V)$. For
$x\in{\bf P}(V)$, the fibre $p^{-1}(x)$  is the linear subspace of ${\cal Q}$
defined by the vanishing of the $Q_\chi$'s, viewed as linear forms in ${\bf
t}$.
These forms are clearly linearly independent as soon as they do not vanish. In
other
words,
 if we denote by
$B_\chi$  the base locus of the quadrics in  $(\sy^2V)_\chi$ and put $B=\cup
B_\chi$, the map $p:{\cal X}\rightarrow {\bf P}(V)$ is a vector bundle
fibration above
${\bf P}(V)\moins B$; in particular ${\cal X}$ is non-singular outside
$p^{-1}(B)$. Therefore it is enough to prove that  ${\cal X}_{\bf t}$ is smooth
at the
points of $B\cap{\cal X}_{\bf t}$.

\ind Observe that an element $x$ in $B$ has two of its
$X$\tx coordinates zero. Since the equations are symmetric in the $X$\tx
coordinates
we may assume $X_\beta=X_\gamma=0$. Then the Jacobian matrix $\displaystyle
({\partial
Q_\chi\over \partial X_\psi}(x))$ takes the form
$$\pmatrix{ 2t^1_1\,X_1&2t^1_2\,X_\alpha&0&0\cr
t^\alpha_1\,X_\alpha & t^\alpha_1\,X_1 &0&0\cr
0&0&t^\beta_1\,X_1&t^\beta_2\,X_\alpha\cr
0&0&t^\gamma_2\,X\alpha&t^\gamma_1\,X_1\cr
}\ .$$
For  generic ${\bf t}$ this matrix is of rank $4$ except when  all the $X$\tx
coordinates of $x$ vanish; but we have seen that this is impossible when ${\bf
t}$ is general enough. \cqfd
\vskip1truecm
{\bf 2. Some comments}
\ind As mentioned in the introduction, the construction is inspired by Reid's
example
of a surface of general type with
$p_g=0$, $K^2=2$, $\pi_1=H$ [R]. This is more than a coincidence. In fact, let
$\widetilde S$ be the hyperplane section $X_1=0$ of $\widetilde X$. It is
stable under
the action of $H$ (so that $H$ acts freely on $\widetilde S$), and one can
prove as
above that it is smooth for a generic choice of the parameters. The surface
$S:=\widetilde S/H$ is a Reid surface,  embedded in $X$ as an ample divisor,
with
$h^0(X,{\cal O}_X(S))=1$. \ind In general, let us consider  a Calabi-Yau
threefold $X$
which contains a {\it rigid ample surface} -- i.e. a smooth ample divisor $S$
such
that  $h^0({\cal O}_X(S))=1$. Put $L:={\cal O}_X(S)$. Then $S$ is a minimal
surface of
general type (because $K_S=L_{|S}$ is ample); by the Lefschetz theorem, the
natural
map  $\pi_1(S)\rightarrow \pi_1(X)$ is an isomorphism. Because of the exact
sequence
$$0\rightarrow {\cal O}_X\longrightarrow L\longrightarrow K_S\rightarrow 0$$the
geometric genus $p_g(S):=h^0(K_S)$ is zero.  \ind One has $K_S^2=L^3$; the
Riemann-Roch theorem on $X$ yields $$1=h^0(L)={L^3\over 6}+{L\cdot c_2\over
12}\ ;$$
by Miyaoka theorem [Mi] one has $L\cdot c_2>0$ (the strict inequality requires
playing around a little bit with the index theorem), hence $K^2_S\le 5$.

 \ind With a few exceptions, the possible fundamental groups of surfaces with
$p_g=0$
and $K_S^2=1$ or $2$ are known (see [B-P-V] for an overview). In the case
$K_S^2=1$,
the fundamental group is cyclic of order $\le 5$; if $K_S^2=2$, it is of order
$\le
9$; moreover the dihedral group $D_8$ cannot occur. I believe that the
symmetric group
${\goth S}_3$ cannot occur either, though I do not think the proof has been
written
down. If this is true,  the quaternionic group $H$ is the only non-Abelian
group
which occurs in this range.
\ind On the other
hand, little is known about surfaces with $p_g=0$ and $K_S^2=3,4$ or 5. Inoue
has
constructed examples with $\pi_1=H\times ({\bf Z}_2)^n$, with $n=K^2-2$ ({\it
loc.~cit.}); I do not know if they can appear as rigid ample surfaces in a
Calabi-Yau
threefold.
\ind Let us denote by $\widetilde X$ the universal cover of $X$,  by
$\widetilde
L$ the pull back of $L$ to $\widetilde X$, and by $\rho$ the representation of
$G$ on
$H^0(\widetilde X,\widetilde L)$. One has $\Tr \rho(g)=0$ for $g\not=1$ by the
holomorphic Lefschetz formula, and $\Tr \rho(1)=\chi(\widetilde
L)=|G|\,\chi(L)=|G|$.
Therefore $\rho$ {\it is isomorphic to the
regular representation}.
 Looking at the list in {\it
loc.~cit.} one gets a few examples of this situation,
 for instance:

 -- $G={\bf Z}_5$, $\widetilde X=\,$a  quintic hypersurface in ${\bf P}^4$;

-- $G=({\bf Z}_2)^3$ or ${\bf Z_4}\times{\bf Z}_2$, $\widetilde X=\,$an
intersection of
$4$ quadrics in ${\bf P}^7$ as above;

-- $G={\bf Z}_3\times{\bf Z}_3$, $\widetilde X=\,$a hypersurface of bidegree
$(3,3)$
in ${\bf P}^2\times{\bf P}^2$.\medskip

\ind Of course when looking for Calabi-Yau threefolds with interesting $\pi_1$
there is
no reason to assume that it contains an ample rigid surface. Observe however
that if
we want to use the preceding method, i.e. find a projective space ${\bf P}(V)$
with an
action of $G$ and a smooth invariant linearly normal Calabi-Yau threefold
$\widetilde
X\i {\bf P}(V)$, then the line bundle ${\cal O}_{\widetilde X}(1)$ will be the
pull-back of an ample line bundle $L$ on $X$, and by the above argument  the
representation of $G$ on $V$ will be  $h^0(L)$ times the regular
representation. This
leaves little hope to find an invariant Calabi-Yau threefold when the product
$h^0(L)\,|G|$ becomes large.
\baselineskip14pt
\vskip 2cm
\def\num#1{\item{\hbox to\parindent{\enskip [#1]\hfill}}}
\def\pc#1{\tenrm#1\sevenrm}
\parindent=1.5cm
\centerline{\bf REFERENCES}
\bigskip
\num{B-P-V} W. {\pc BARTH}, C. {\pc PETERS}, A. {\pc VAN DE} {\pc VEN}: {\sl
Compact
complex surfaces}. Ergebnisse der Math., Springer-Verlag (1984). \smallskip
\num{Mi} Y. {\pc MIYAOKA}: {\sl The Chern classes and Kodaira dimension of a
minimal
variety}. Adv. Studies in Pure Math. {\bf 10}, 449-476, Kinokuniya-North
Holland
(1987). \smallskip
\num{R} M. {\pc REID}: {\sl Surfaces with $p_g=0$, $K^2=2$}. Unpublished
manuscript (1979).\vskip1cm \hfill\hbox to 5cm{\hfill A. Beauville\hfill}\par
\hfill\hbox to 5cm{\hfill URA 752 du CNRS\hfill}\par
\hfill\hbox to 5cm{\hfill Math\'ematiques -- B\^at. 425\hfill}\par
\hfill\hbox to 5cm{\hfill Universit\'e Paris-Sud\hfill}\par
\hfill\hbox to 5cm{\hfill 91 405 {\pc ORSAY} Cedex, France\hfill}\par

\bye